\documentclass[12pt]{amsart}
\usepackage{amssymb}
\usepackage{graphicx}

\newtheorem{theorem}{Theorem}[section]

\title[Skolem's conjecture for a family of exponential equations]{Skolem's conjecture confirmed for a family of exponential equations, II$^1$\footnote{\lowercase{1 \uppercase{T}he paper appeared in \uppercase{A}cta \uppercase{A}rith. 197.2 (2021), 129--136. \uppercase{H}ere we correct a typo: on page 134 in line 5, the modulus $m'$ should read as $m'=|y^{s_N}(x^{o_N}-1)|$. (\uppercase{L}ike here on page 6 in line 7 up.) \uppercase{W}e are very much grateful to \uppercase{R}eese \uppercase{S}cott and \uppercase{R}obert \uppercase{S}tyer for pointing out this error.}}}

\author[A. B\'erczes]{A. B\'erczes}
\author[L. Hajdu]{L. Hajdu}
\author[R. Tijdeman]{R. Tijdeman}

\thanks{Research was supported in part by grants 115479, 128088 and 130909 of the Hungarian National Foundation for Scientific Research and by the projects EFOP-3.6.1-16-2016-00022 and EFOP-3.6.2-16-2017-00015, co-financed by the European Union and the European Social Fund.}



\address{A. B\'erczes, L. Hajdu \newline
         \indent Institute of Mathematics\newline
         \indent University of Debrecen \newline
         \indent H-4010 Debrecen, P.O. Box 12, Hungary}
\email{berczesa@science.unideb.hu}
\email{hajdul@science.unideb.hu}

\address{R. Tijdeman \newline
         \indent Mathematical Institute\newline
         \indent Leiden University\newline
         \indent Postbus 9512, 2300 RA Leiden, The Netherlands}
\email{tijdeman@math.leidenuniv.nl}

\begin{document}

\dedicatory{Dedicated to K\'alm\'an Gy\H{o}ry on the occasion of his 80th birthday.}

\begin{abstract} According to Skolem's conjecture, if an exponential Diophantine equation is not solvable, then it is not solvable modulo an appropriately chosen modulus. Besides several concrete equations, the conjecture has only been proved for rather special cases. In this paper we prove the conjecture for equations of the form $x^n-by_1^{k_1}\dots y_\ell^{k_\ell}=\pm 1$, where $b,x,y_1,\dots,y_\ell$ are fixed integers and $n,k_1,\dots,k_\ell$ are non-negative integral unknowns. This result extends a recent theorem of Hajdu and Tijdeman.
\end{abstract}

\maketitle

\section{Introduction}

Skolem \cite{sk} proposed the following conjecture: if a purely exponential Diophantine equation is not solvable, then it is not solvable modulo an appropriate modulus. Skolem's conjecture and variants of it (besides several concrete cases) have been proved only under restricted conditions. Schinzel \cite{sc1} (extending results of Skolem \cite{sk}) verified the conjecture in the one-term case, that is for equations of the shape $\alpha_1^{k_1}\cdots \alpha_\ell^{k_\ell}=\beta$. Here $\alpha_1,\dots,\alpha_\ell$ and $\beta$ are fixed elements of a number field, and $k_1,\dots,k_\ell$ are unknown integers. Bartolome, Bilu and Luca \cite{bbl} proved Skolem's conjecture for equations of the form $\lambda_1 \alpha_1^k+\dots+\lambda_\ell \alpha_\ell^k=0$, where $\lambda_1,\dots,\lambda_\ell$ and $\alpha_1,\dots,\alpha_\ell$ are elements of a fixed number field such that the order of the multiplicative group generated by $\alpha_1,\dots,\alpha_\ell$ is one and $k$ is an integral variable. The main result of \cite{bbl} follows from those in \cite{sc1,sc2} - though some nontrivial argument is needed \cite{sc4}. Bert\'ok and Hajdu \cite{bh1,bh2} proved that in some sense Skolem's conjecture is valid for 'almost all' equations. For strongly related problems and results concerning recurrence sequences, see the papers \cite{sc2,sc3,os}, and the references there.

Recently, Hajdu and Tijdeman \cite{ht} proved that Skolem's conjecture is valid for the Catalan equation $x^n-y^k=1$, if $x$ and $y$ are fixed positive integers one of which is prime and $n,k$ are non-negative integer variables. Further related theoretical and numerical results can be found e.g. in \cite{sc2,bl,be,bh1,bh2}, see also the references there. In the present paper we prove Skolem's conjecture for equations of the form $x^n-by_1^{k_1}\cdots y_\ell^{k_\ell}=\pm 1$, where $b,x,y_1,\dots,y_\ell$ are fixed integers and $n,k_1,\dots,k_\ell$ are non-negative integral unknowns. This result extends the above mentioned theorem of Hajdu and Tijdeman \cite{ht} that covers the case $b=\ell=1$ and one of $x,y_1$ is a prime. Note that in particular, our present results imply that Skolem's conjecture holds for the Catalan equation $x^n-y^k=1$ for fixed bases $x,y$. The same is valid for another well known equation, namely for $\frac{x^n-1}{x-1}=y^k$.

\section{The theorem and its proof}

Our main result is the following. In its formulation, as well as in the rest of the paper, we use the convention $0^0=0$.

\begin{theorem}
\label{thm1}
Let $b,x,y_1,\dots,y_\ell$ be integers. Then there exists a modulus $m$ such that the congruence
\begin{equation}
\label{maincong}
x^n-by_1^{k_1}\cdots y_\ell^{k_\ell}\equiv \pm 1\pmod{m}
\end{equation}
has precisely the same solutions in non-negative integers $n,k_1,\dots,k_\ell$ as the equation
\begin{equation}
\label{maineq}
x^n-by_1^{k_1}\cdots y_\ell^{k_\ell}=\pm 1
\end{equation}
has.
\end{theorem}

\noindent{\bf Remark 1.} In the statement of Theorem \ref{thm1} the $b$'s can be omitted. They are included for its proof.

\vskip.2cm

\noindent{\bf Remark 2.} It will be clear from the proof that given $b, x,y_1,\dots,y_\ell$, the modulus $m$ can be explicitly constructed, and can be bounded in terms of $b, x,y_1,\dots,y_\ell$.

\vskip.2cm

\noindent{\bf Remark 3.}
For every $m$ all solutions of \eqref{maineq} are solutions of \eqref{maincong}. So it suffices to prove that for certain $m$ every solution of \eqref{maincong} is a solution of \eqref{maineq}. To find such an $m$ the following observation will be useful. For fixed $b,x,y_1,\dots,y_\ell$, write $S_\infty$ for the set of solutions of \eqref{maineq} and for any modulus $m$ let $S_m$ be the set of solutions of \eqref{maincong}. Then we have $S_\infty\subseteq S_m$ for any $m\geq 2$. On the other hand, if $m_1,m_2\mid m$ then we clearly have $S_m\subseteq S_{m_1}\cap S_{m_2}$. These altogether imply that if we can find moduli $m_1,\dots,m_t$ such that
$$
\bigcap\limits_{i=1}^t S_i=S_\infty
$$
then the choice
$$
m:=\prod\limits_{i=1}^t m_i
$$
is appropriate: we have $S_\infty=S_m$, that is, for this $m$ the solutions of \eqref{maincong} and \eqref{maineq} coincide. Furthermore, if we have found a modulus $m'$ such that the terms of all the solutions of \eqref{maincong} with modulus $m'$ are bounded (either because the bases are $\pm 1$ or $0$, or because the unknown exponents are bounded), then we may choose $m''$ sufficiently large so that for the modulus $m=m'm''$ \eqref{maincong} and \eqref{maineq} have exactly the same solutions. To see this, observe that if
$$
A\equiv B\pmod{m''}
$$
holds with some $A,B\in{\mathbb Z}$ such that $|A|+|B|<m''$, then clearly $A=B$. 

This is the strategy we follow in the proof. First we find moduli $m_i$ such that $S_{m'}$ is bounded for the product $m'$ of these moduli. (Note that taking $m'$ to be the lcm of the $m_i$-s is an appropriate choice, too.) Then we take an $m''$ large enough and choose $m=m'm''$. Note that if, in particular, with some $m$ we get that $S_m$ is empty, then so is $S_\infty$, and our statement immediately follows.

We illustrate our method by the following simple example. Consider the equation
\begin{equation}
\label{ujeq}
5^x-2^y=3
\end{equation}
and the corresponding congruence
\begin{equation}
\label{ujcong}
5^x-2^y\equiv 3\pmod{m}.
\end{equation}
Then the sets of solutions of \eqref{ujcong} with $m=4$ and $m=25$ are
$$
S_4=\{(u,1):u\geq 0\}
$$
and
$$
S_{25}=\{(0,20v+11),\ (1,20v+1),\ (u,20v+17):u\geq 2,\ v\geq 0\},
$$
respectively. Hence
$$
\{(1,1)\}\subseteq S_\infty\subseteq S_{100}=S_{4}\cap S_{25}=\{(1,1)\}.
$$
Therefore $S_{100}=S_{\infty}=\{(1,1)\}$.

\vskip.2cm

\noindent{\bf Remark 4.}
Our theorem implies that for fixed $x,y$ all solutions of the Catalan equation
$$
x^n-y^k=1
$$
can be found locally. Hence, for any positive integers $x,y>1$ with $(x,y)\neq (3,2)$ one can find a modulus $m=m(x,y)$ such that the congruence
$$
x^n-y^k\equiv 1\pmod{m}
$$
has no solutions in integers $n,k>1$.

A similar remark applies for another well known equation, namely for
\begin{equation}
\label{lasteq}
\frac{x^n-1}{x-1}=y^k
\end{equation}
in integers $x,y,n,k$ with $x>1$, $y>1$, $n>2$, $k>1$. For the history of \eqref{lasteq} and for related results see e.g. \cite{st0} or Theorem 12.5 of \cite{st}. Obviously, our theorem implies that for any fixed $x,y$ there exists an $m$ such that \eqref{lasteq} and \eqref{lasteq} modulo $m$ have the same solutions.

\vskip.2cm

\begin{proof}[Proof of Theorem \ref{thm1}]
We start with some trivial cases.

Assume first that $x=0$. If we also have $by_1\dots y_\ell=0$, then \eqref{maincong} is not solvable with $m=2$ and we are done. Otherwise, let $p_0$ be a prime factor of $b$, or $p_0=1$ if $|b|=1$, and for $i=1,\dots,\ell$ let $p_i$ be a prime factor of $y_i$ with $p_i=1$ whenever $|y_i|=1$. Put $m'=p_0p_1\cdots p_\ell$. If $m'=1$ then the numbers $b,y_1,\dots,y_\ell$ are all $\pm 1$, and all solutions of congruence \eqref{maincong} with $m=3$ are solutions of equation \eqref{maineq} itself. If $m'>1$ then \eqref{maincong} with $m=m'$ has only solutions where the exponents of the bases $y_i$ with $|y_i|>1$ are zero. Thus \eqref{maincong} has no solutions with $m=3m'$ in this case, and we are done. From now on we assume $x\neq 0$.

Suppose next that $by_1\dots y_\ell=0$. If $|x|=1$, then equation \eqref{maineq} has the same solutions as \eqref{maincong} with $m=3$. If $|x|>1$ then \eqref{maincong} with $m=|x|$ implies $n=0$. It follows that \eqref{maineq} and \eqref{maincong} with $m=3|x|$ have the same solutions.

From this point on we suppose $bxy_1\cdots y_\ell\neq 0$. We apply induction on $\ell$. Consider first $\ell=0$, i.e. $y_1^{k_1}\cdots y_\ell^{k_\ell}=1$. If $-b$ equals the right hand side of \eqref{maineq}, then we choose $m=p$ where $p$ is a prime not dividing $x$ to see that \eqref{maincong} has no solutions. So assume that $-b$ is different from the right hand side of \eqref{maineq}. Then if $|x|=1$, we choose $m=|b|+3$ and check that \eqref{maincong} and \eqref{maineq} have the same solutions. If $|x|>1$, we take $m'=|x|^{n_0}$ where $n_0$ is an integer with $|x|^{n_0}>2(|b|+1)$. Then we see that \eqref{maincong} with $m=m'$ has no solutions for $n\geq n_0$. Hence (as we explained in Remark 3) we can easily find an $m''$ such that \eqref{maineq} and \eqref{maincong} with $m=m'm''$ have the same solutions.  

Let $\ell=1$. For simplicity, we rewrite equation \eqref{maineq} as
\begin{equation}
\label{maineq2}
x^n-by^k=c
\end{equation}
with $c=\pm 1$. If $|y|=1$, we can adapt the above argument for $\ell=0$. So we may assume $|y|>1$. If $|x|=1$, then $0<|by^k|=|x^n-c|=2$, and modulo $|y|^2$ we see that $k\leq 1$. But then both $x^n$ and $by^k$ are bounded, and we can easily find a modulus $m$ such that \eqref{maincong} and \eqref{maineq} have the same solutions. If $\gcd(x,by)=d>1$ then \eqref{maincong} with $m=d$ gives that $n=0$ or $k=0$. However, these cases are equivalent with the already treated cases $x=1$ and $y=1$, respectively. Altogether, in the sequel without loss of generality we may assume that $|x|>1$, $|y|>1$ and $\gcd(x,by)=1$.

The general argument for $\ell=1$ needs some preparation. Put
$$
S=\{p:\ p\ \text{is a prime divisor of}\ bxy\}
$$
and write $U_S$ for the set of integers having all their prime divisors in $S$. Consider the equation
\begin{equation}
\label{sueq}
v_1-v_2=c.
\end{equation}
It is well-known that \eqref{sueq} has only finitely many solutions in $v_1,v_2\in U_S$ whose number can be effectively bounded in terms of $b,c,x,y$. For such theorems see e.g. the book of Evertse and Gy\H{o}ry \cite{egy}. Hence the number of solutions of
\begin{equation}
\label{sueq2}
x^t-z=c
\end{equation}
in integers $t,z$ such that $t\geq 0$ and $z\in U_S$ is also bounded in terms of $b,c,x,y$. 
Write
$$
N=N(b,c,x,y)
$$
for the number of values $t$ appearing in the solutions of equation \eqref{sueq2}.

Now we turn to our main argument. At this point we distinguish two cases. Assume first that $c=1$ in \eqref{maineq2}. Let $s_1$ be the smallest integer such that
$$
|y|^{s_1}>|x|+1.
$$
Obviously, $s_1$ can be easily expressed in terms of $x,y$. 
If $k<s_1$, then $k$ is bounded and can be considered to be fixed. By that $by^k$ will be fixed and we are back to the case $\ell = 0$. So we may suppose $k\geq s_1$. Then we get
$$
x^n\equiv 1\pmod{|y|^{s_1}}.
$$
Thus the order $o_1$ of $x$ modulo $|y|^{s_1}$ must divide $n$. By our choice of $s_1$ we have that this order is not one, so
$$
2\leq o_1\leq |y|^{s_1}.
$$
Let now $s_2$ be the smallest integer such that
$$
|y|^{s_2}>|x|^{o_1}+1.
$$
Observe that $o_1$ and $s_2$ are bounded in terms of $x,y$. If $k<s_2$ we can proceed as in the case $k<s_1$. So we may assume that $k\geq s_2$. Hence we obtain
$$
x^{o_2}\equiv 1\pmod{|y|^{s_2}}.
$$
Therefore the order $o_2$ of $x$ modulo $|y|^{s_2}$ must also divide $n$. We have $o_1\mid o_2$, too. Further, by our choice of $s_2$ we see that
$$
1<o_1<o_2\leq |y|^{s_2}.
$$
Continuing this procedure, we have two options. Either the process terminates in at most $N$ steps, yielding modulo $|y|^{s_i}$ for some $i\leq N$ that $k$ is bounded in terms of $b,c,x,y$. Then we are back in the case $\ell=0$ and we are done. Or, after $N$ steps we obtain that there exist divisors $o_1,\dots,o_N$ of $n$ with
$$
1<o_1<\dots<o_N\leq |y|^{s_N}
$$
where $s_N$ is bounded in terms of $x,y$, such that 
$$
o_1\mid o_2,\hdots, o_{N-1}\mid o_N,o_N\mid n.
$$
Put $o_0=1$ and consider \eqref{maineq2} modulo $x^{o_i}-1$ for $i=0,1,\dots,N$. We get that for $k\geq {s_N}$
$$
by^k\equiv 0\pmod{x^{o_i}-1}
$$
holds, hence $x^{o_i}-1\in U_S$ $(i=0,1,\dots,N)$. It follows that
$$
x^{o_i}-z_i=1
$$
holds for $i=0,1,\dots,N$ for some $z_i\in U_S$. However, this contradicts the definition of $N$. So taking the modulus
$$
m'=|y^{s_N}(x^{o_N}-1)|
$$
we get that in all solutions of \eqref{maineq2} modulo $m'$, we have $k<{s_N}$. Then, as we already mentioned, we are back to the case $\ell=0$, and our claim follows for $c=1$.

Assume now that $c=-1$ in \eqref{maineq2}. If $by$ has no odd prime factors, then our equation reduces to an equation of the shape
$$
x^n\pm 2^v=-1
$$
in integers $n,v$ with $n\geq 0$ and $v\geq v_0$, where $v_0$ is the non-negative integer satisfying $b=\pm 2^{v_0}$. This we can rewrite as
$$
2^v\pm x^n=\mp 1.
$$
If we have $1$ on the right hand side, then we are done by our results for $c=1$. On the other hand, if $-1$ stands on the right hand side, then by $|x|>1$ and $\gcd(x,by)=1$ we know that $x$ has an odd prime factor. This shows that by interchanging the roles of $x$ and $y$ if necessary, without loss of generality we may assume that $by$ has an odd prime factor $p$. Write $p-1=2^r h$ with $r\geq 0$ and $h$ odd. Considering \eqref{maineq2} modulo $p$, we get that either $k=0$, whence we are back to the already treated case $\ell=0$, or
$$
x^n\equiv -1\pmod{p}.
$$
This implies that writing $n=2^s n'$ with $s\geq 0$ and $n'$ odd, $s<r$ holds. Indeed, otherwise by the above congruence and $p\nmid x$ we would get
$$
1\equiv \left(x^{p-1}\right)^{2^{s-r} n'}\equiv \left(x^{2^r h}\right)^{2^{s-r} n'}\equiv \left(x^{2^s n'}\right)^h\equiv -1\pmod{p}, 
$$
which is a contradiction. This means that we can rewrite \eqref{maineq2} as
\begin{equation}
\label{eqq}
x_s^{n'}-by^k=-1,
\end{equation}
where $x_s=x^{2^s}$ for some $s$ with $0\leq s<r$. Clearly, $x_s$ can be considered to be fixed. Note that the exponent $n'$ is odd. Now considering the above equation modulo appropriate powers of $|y|$, since $\text{ord}_{|y|^i}(x_s)\mid 2n'$ for every $i\geq 1$, similarly as in case of $c=1$, we obtain divisors $o_i^*$ of $2n'$, hence odd divisors $o_i'$ of $n'$, such that $1=o_0'<o_1'<o_2'<\dots$ and also $o_0'\mid o_1',o_1'\mid o_2',o_2'\mid o_3',\dots$. As for all such $o_i'$ we have
$$
x^{{o_i'}}+1\mid x^{n'}+1,
$$
considering equation \eqref{eqq} modulo $x^{o'_N}+1$, we get too many solutions $t$ to the $S$-unit equation
$$
x^t-z=-1,
$$
unless $k$ is bounded modulo some fixed power of $|y|$, in terms of $b,c,x,y$. Hence we are back to the situation $\ell=0$, and our claim follows also in this case. Hence the proof concerning the case $\ell=1$ is complete.

Let now $\ell\geq 2$ be arbitrary, and suppose that the statement holds for all $\ell'$ with $\ell'<\ell$ and for all $b,x,y_1,\dots,y_{\ell'}$. We keep the notation from the case $\ell=1$ with a simple change: $S$ is replaced by
$$
S=\{p:\ p\ \text{is a prime divisor of}\ bxy_1\cdots y_\ell\},
$$
and $U_S$ and $N$ are also re-defined with this choice of $S$. 

We assume that we have $1$ on the right hand sides of \eqref{maincong} and \eqref{maineq}; the treatment of the other case follows similar lines, with certain small modifications as for $\ell=1$. We use a similar construction as in case of $\ell=1$. First, let $r_1$ be the smallest integer such that
$$
|y_1\cdots y_\ell|^{r_1}>|x|+1.
$$
Clearly, $r_1$ can be bounded in terms of $x,y_1,\dots,y_\ell$. If there exists an $i\in\{1,\dots,\ell\}$ such that
$k_i<r_1$ then $y_i^{k_i}$ can be merged into $b$, and we are done by the induction hypothesis. So we may assume that $r_1\leq k_i$ for all $i=1,\dots,\ell$. Then modulo $|y_1\cdots y_\ell|^{r_1}$ equation \eqref{maineq} yields that the order $w_1$ of $x$ modulo $|y_1\cdots y_\ell|^{r_1}$ divides $n$. Repeating this argument $N$ times, just as in the proof of the case $\ell=1$, we obtain integers $1=w_0<w_1<\dots<w_N$ with
$$
w_0\mid w_1,w_1\mid w_2,\dots,w_{N-1}\mid w_N,w_N\mid n,
$$
all bounded above in terms of $x,y_1,\dots,y_\ell$. Then the modulus
$$
|y_1\cdots y_\ell|^{r_N}
$$
provides an upper bound for $\max(k_1,\dots,k_\ell)$. Hence by the induction hypothesis our claim follows also in this case. This completes the proof of our statement.
\end{proof}

\section*{Acknowledgements}

The authors are grateful to the Referee for the insightful and helpful remarks.

\end{document}